\documentclass[12pt]{article}
\usepackage{amsxtra}
\usepackage{amssymb, amsmath}
\usepackage{amsthm}
\usepackage {graphicx}

\begin{document}
\begin{center}
\textbf{\LARGE{Universal Bimagic Squares and \\ the day 10}$^{th}$\textbf{ October 2010
(10.10.10)}}
\end{center}

\bigskip
\begin{center}
\textbf{\large{Inder Jeet Taneja}}\\
Departamento de Matem\'{a}tica\\
Universidade Federal de Santa Catarina\\
88.040-900 Florian\'{o}polis, SC, Brazil.\\
\textit{e-mail: ijtaneja@gmail.com\\
http://www.mtm.ufsc.br/$\sim$taneja}
\end{center}

\bigskip
\begin{abstract}

\textit{In this short note we have produced for the first time in the history different kinds of universal bimagic squares. This we have made using only the digits 0,1 and 2. The universal bimagic squares of order $8 \times 8$ and $16 \times 16$ are with the digits 0 and 1. The universal bimagic square of order $9 \times 9$ is with the digits 0, 1 and 2. It is interesting to note that the day 0ctober 10 have only the digits 0 and 1 if we consider it as 10.10.10. If we consider the date as 10.10.2010, then this has the digits 0, 1 and 2.}
\end{abstract}

\section{Details}

In this work we shall present \textit{universal bimagic} squares of order $8 \times 8$ and
$16 \times 16$ having only the digits 0 and 1. A \textit{universal bimagic} square of order
$9 \times 9$ is also presented having three digits 0, 1 and 2. These magic squares are based on
the date October 10 (10.10.10 or 10.10.2010).

\subsection{Universal and Bimagic Squares}

Here below are some definitions.

\bigskip
\noindent \textbf{$\bullet$ Magic square}

\smallskip
\noindent  A magic square is a collection of numbers put as a square matrix,
where the sum of elements of each row, sum of elements of each column or sum
of elements of each of two principal diagonals are the same. For
simplicity, let us write this sum as \textbf{S1}.

\bigskip
\noindent \textbf{$\bullet$ Bimagic square}

\smallskip
\noindent Bimagic square is a magic square where the sum of square of each
element of rows, columns or two principal diagonals are the same. For
simplicity, let us write this sum as \textbf{S2}.
\bigskip

\noindent \textbf{$\bullet$ Universal magic square}

\smallskip
\noindent Universal magic square is a magic square with the following properties:

\begin{itemize}
\item[(i)] \textbf{Upside down}, i.e., if we rotate it to 180 degrees, it remains magic square again;
\item[(ii)] \textbf{Mirror looking}, i.e., if we put it in front of mirror or see from
the other side of the glass, or see on the other side of the paper, it
always remains the magic square.
\end{itemize}

\subsection{The date 10.10.10}

It is interesting to note that at 10 hours, 10 minutes and 10 seconds of the
day 10, month 10 and the year 10 have only the digits 1 and 0, i.e.,
10-10-10-10-10-10. Let us divide it in two parts, i.e., 101010 -- 101010.
Thus we have two equal blocks of six algarisms. If we go only for the day,
these digits repeats on others days too, such as

\begin{center}
01-10-10; 01-01-10; 10-01-10, 11-10-10, etc.
\end{center}

If we go on hours, minutes and seconds we have many combinations of six
algarisms only with the digits 0 and 1.

\bigskip

\noindent \textbf {\large{$\bullet$ $8 \times 8$ -- Universal bimagic square of binary digits 0 and 1}}
\bigskip

\noindent We can make $2\times 2\times 2\times 2\times 2\times 2$ or $2^6 = 64$ different
numbers of six algarisms with the digits 0 and 1. Also, we can write $64 = 8\times 8$.
Here below is a \textbf{universal bimagic square} of
order $8 \times 8$ having 64 different numbers using only the digits 0 and 1.

\begin{figure}[htbp]
\centerline{\includegraphics[width=6.00in,height=3.00in]{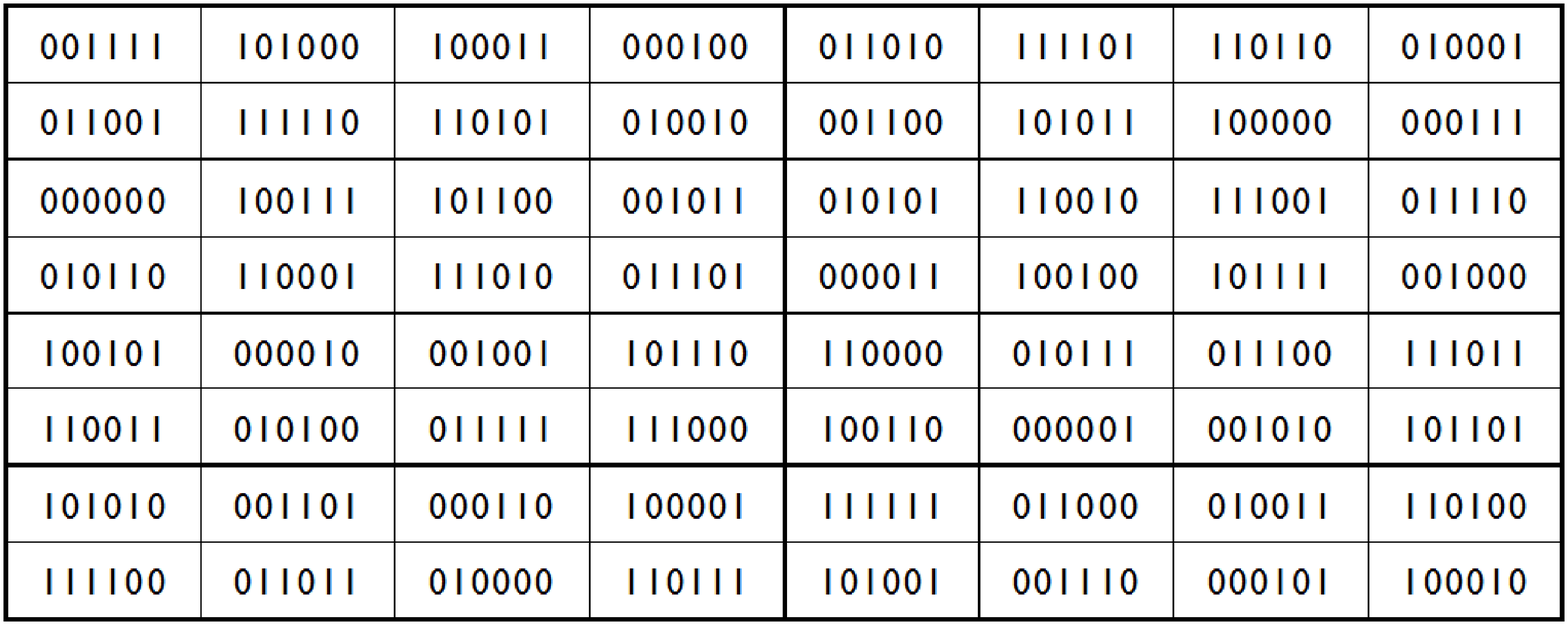}}
\label{fig1}
\end{figure}

\noindent
S1:=44444\\
S2:=44893328844

\bigskip
Also we have sum of each block of $2\times 4=44444$ and the square of sum of each
term in of each block of $2 \times 4=44893328844$

\bigskip
\noindent \textbf {\large {$\bullet$ $16 \times 16$ -- Universal bimagic square of binary digits 0 and 1}}

\bigskip
\noindent Instead, considering six algarisms using only the digits 0 and 1, if we
consider eight algarisms using we can make $2^8 = 256$ different numbers only with the digits 0 and 1. Also we we can write 256 as $16 \times 16$. Here below is a universal bimagic square of order $16 \times 16$ with these 256 different
numbers made from the digits 0 and 1:

\begin{figure}[htbp]
\centerline{\includegraphics[width=6.20in,height=4.00in]{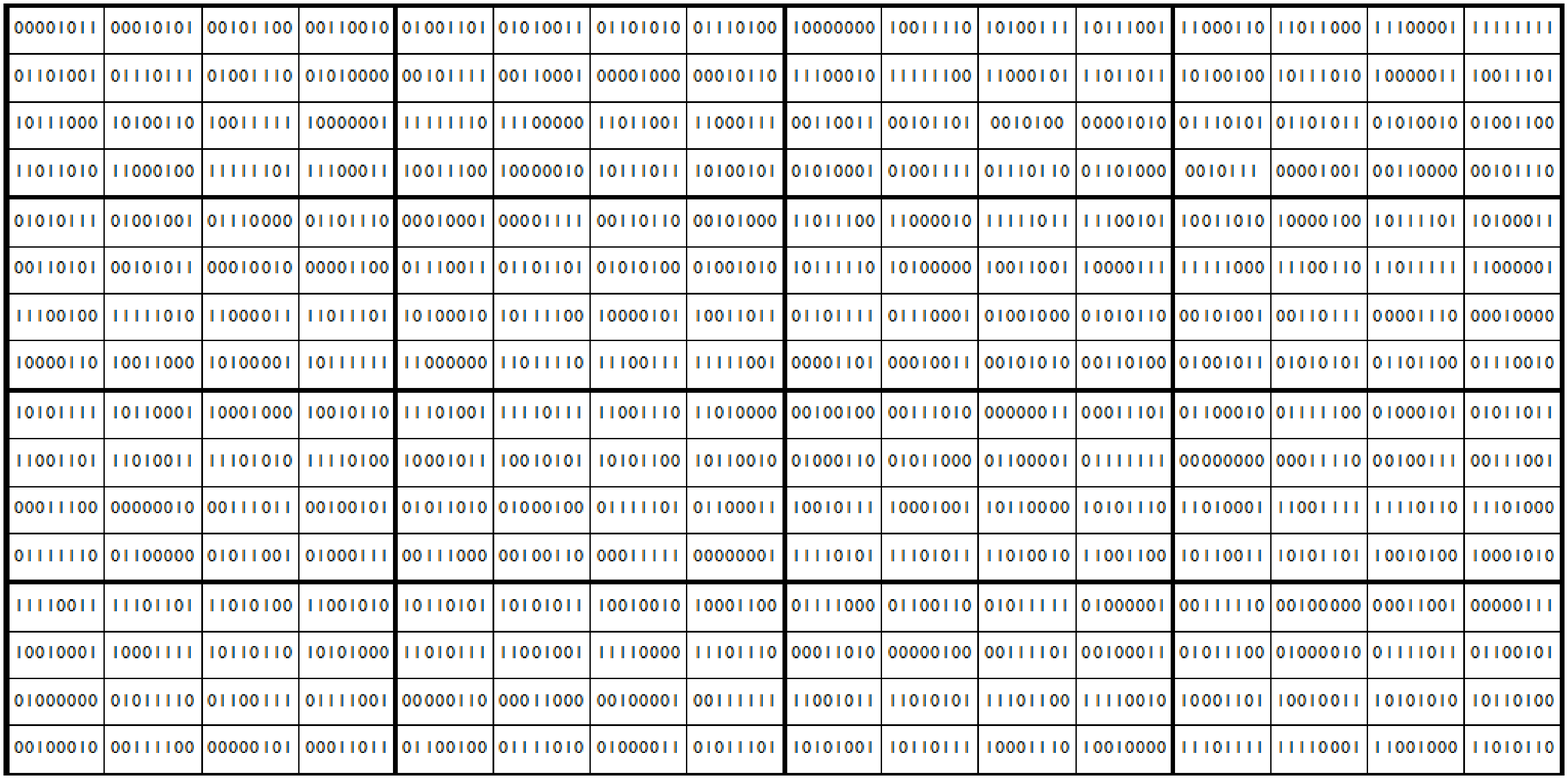}}
\label{fig2}
\end{figure}

\noindent
S1:=88888888\\
S2:=897867554657688

\bigskip
Also we have sum of each block of $4 \times 4=88888888$ and square of sum of each term
in of each block of $4 \times 4=897867554657688$.

\subsection{The date 10.10.2010}

Instead,  considering the year as 10, if we consider it as 2010, then we have
three algarisms 0, 1 and 2. These digits happens on other days too, such as

\begin{center}
02.01.2010; 02.02.2010; 20.10.2010; 02.10.2010; 12.10.2010;  2.10.2010, etc.
\end{center}

Still there are many other dates having only the digits 0, 1 and 2.

\newpage

\noindent \textbf{\large{$\bullet$ $9 \times 9$  -- Universal bimagic square of digits 0, 1 and 2}}

\bigskip

\noindent We can make exactly 81 different numbers having four algarisms from the three digits 0, 1 and 2, i.e, $3\times 3\times 3\times 3 = 81$. Also we can write, $81 = 9\times 9$. Here below is a universal bimagic square of order $9 \times 9$ having only the digits 0, 1 and 2 with 81 different numbers.

\bigskip

\begin{figure}[htbp]
\centerline{\includegraphics[width=6.00in,height=3.00in]{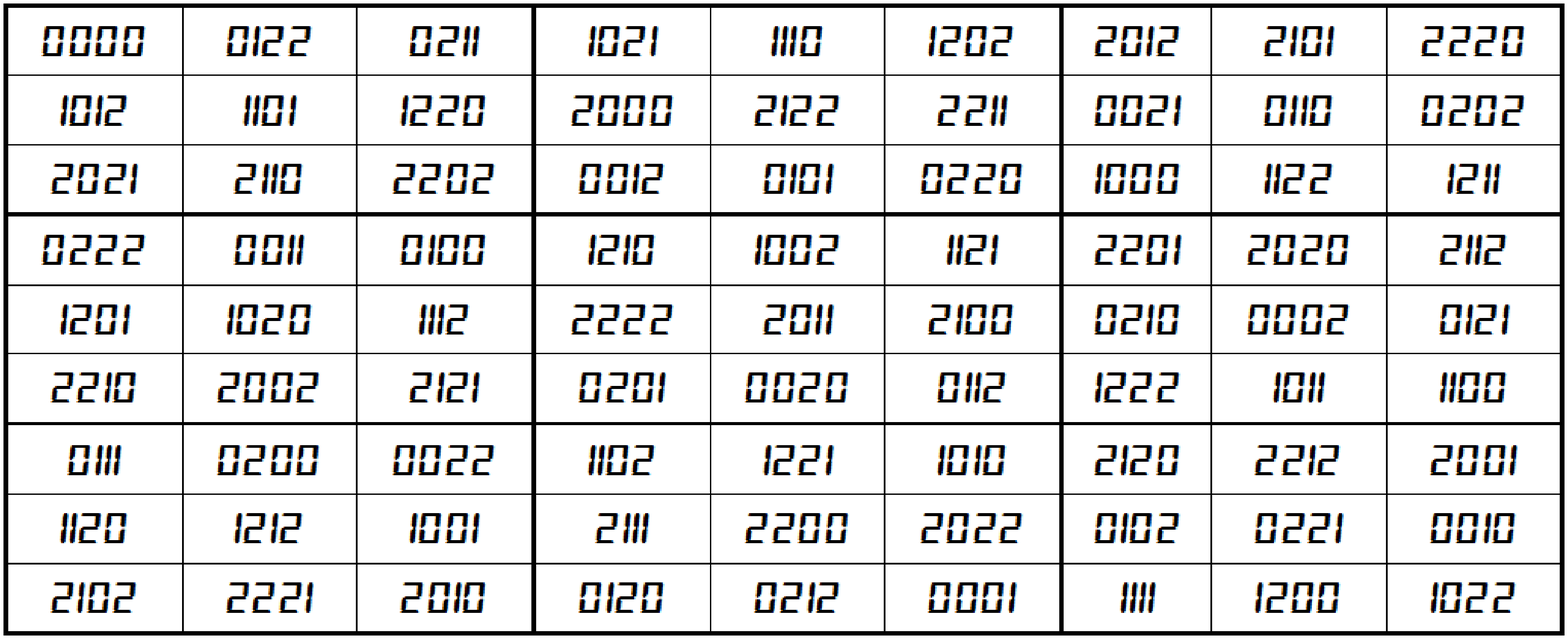}}
\label{fig3}
\end{figure}

\noindent
S1:=9999\\
S2:=17169395

\bigskip
Also we have sum of each block of $3\times 3=9999$ and square of sum of each term in
of each block of $3 \times 3=17169495$

\bigskip
We observe that from the above magic square that if we make a rotation of 180
degrees the digits 2 remains the 2 but if we see it in the mirror 2 becomes
5. Obviously, in this case the sum S1 and S2 are not the same as given above. But still
it is a magic square. If we want to have the same sum, we have to use
2 and 5 together (in the digital form) with either 0, 1 or 8. This study is given in the anther
work Taneja \cite{tan3}.

\bigskip
For more studies on magic and bimagic squares, we suggest to the readers
the two sites \cite{boyer}, \cite{her} where one can find a good collection of work,
papers, books, etc. The idea of universal bimagic square is presented for the first time here.

\begin{center}
---------------------------
\end{center}

\end{document}